\documentclass[11pt,a4paper,twoside,draft]{amsart}

\usepackage{amsmath}
\usepackage{amsfonts}
\usepackage{amssymb}

\newcommand{\beqno}{\begin{eqnarray*}}
\newcommand{\eeqno}{\end{eqnarray*}}

\newcommand{\mint}{-\hskip-12pt\int}
\newcommand{\smint}{-\hskip-10.3pt\int}

\newcommand{\beqla}[1] {\begin {eqnarray}\label{#1}}
\def \eeq {\end {eqnarray}}
\newtheorem{theorem}{Theorem}
\newtheorem{lemma}[theorem]{Lemma}
\newtheorem{proposition}[theorem]{Proposition}
\newtheorem{corollary}[theorem]{Corollary}

\theoremstyle{definition}
\newtheorem{definition}[theorem]{Definition}
\newtheorem{example}[theorem]{Example}

\theoremstyle{remark}
\newtheorem{remark}[theorem]{Remark}

\newcommand{\boundary}{\partial}

\newcommand{\real}{{\mathbf R}}

\newcommand{\integer}{{\mathbf Z}}
\newcommand{\nanu}{{\mathbf N}}
\newcommand{\sdist}{{\mathcal S'}}
\newcommand{\distr}{{\mathcal D'}}
\newcommand{\ssmooth}{{\mathcal S}}
\newcommand{\hlm}{{{\mathcal M}}}
\newcommand{\ds}{\displaystyle}
\newcommand{\smoothie}[1]{C^\infty_{0,#1}}
\newcommand{\capa}{{\rm cap}}

\newcommand{\refeq}[1]{(\ref{#1})}

\DeclareMathOperator{\supp}{supp}
\newcommand{\smooth}{C_0^\infty}

\DeclareMathOperator{\diam}{{\rm diam\,}}

\begin{document}

\title[characterizations of Hardy-Sobolev functions]
{Pointwise characterizations of Hardy-Sobolev functions}

\author{Pekka Koskela}
\address{University of Jyv\"askyl\"a, Department of Mathematics and Statistics,
P.O.~Box~35 (MaD), FIN-40014 University of Jyv\"askyl\"a,
Finland}
\email{pkoskela@maths.jyu.fi}

\author{Eero Saksman}
\address{University of Jyv\"askyl\"a, Department of Mathematics and Statistics,
P.O.~Box~35 (MaD), FIN-40014 University of Jyv\"askyl\"a,
Finland}
\email{saksman@maths.jyu.fi}





\begin{abstract} We establish  pointwise characterizations of functions in the  Hardy-Sobolev spaces $ H^{1,p}$ 
within the range $p\in (n/(n+1),1]$. In particular, a locally integrable function $u$ belongs
to $ H^{1,p}(\real^n)$ if and only if  $u\in L^p(\real^n)$ and it satisfies the Hajlasz type condition
$$
|u(x)-u(y)|\leq |x-y|(h(x)+h(y)),\quad x,y\in \real^n\setminus E,
$$
where $E$ is a set of measure zero and $h\in L^p(\real^n)$.
We also investigate  Hardy-Sobolev spaces on subdomains and  extend Hardy inequalities to the case $p\leq 1.$ 
\end{abstract}

\maketitle

\section{Introduction}\label{se:intro}

It is a well-established fact that, for the purposes of harmonic analysis or 
theory of partial differential equations, the right substitute for $L^p(\real^n)$ in case $p\in (0,1]$ is
the (real) Hardy space $H^p(\real^n)$, or it's local version $h^p(\real^n)$. The Hardy spaces, 
or  their local versions if needed, behave nicely
under the action of regular singular integrals or  pseudo-differential operators. Moreover, in the case of
 Hardy spaces the Paley-Littlewood
theory and interpolation results extend to the whole scale of Lebesgue exponents $p\in (0,\infty).$ It is
hence natural to investigate Sobolev spaces where one (roughly speaking) demands that the $s$:th derivative
belongs to a Hardy type space in the case $p\leq 1$. After the fundamental work of Fefferman and Stein 
\cite{FeffermanStein} this line
of research 
was initiated by Peetre in early 70's,  and it was  generalized and carried further by Triebel and others.
We refer to \cite{Peetre}, \cite{Triebel} for extensive accounts on  general Besov and Triebel -type
scales of function spaces in the case $p\in (0,1].$

In this paper we establish new pointwise characterizations of Hardy-Sobolev spaces in the most
important case where the smoothness index is one and the elements in these spaces are honest functions, i.e. they belong to
$L^{1}_{loc}$.
Recall  that a distribution $f\in S'(\real^n)$ belongs to the homogeneous 
(Hardy-)Sobolev space $\dot H^{1,p}(\real^n)$ if $D_kf\in H^p(\real^n)$ for
$1\leq k\leq n.$ Modulo polynomials, these spaces coincide with the homogeneous spaces considered in \cite[Chapter 5]{Triebel}.
Various characterizations in terms of Paley-Littlewood decomposition
(the square function), Lusin functions, atoms, maximal operators, or various integral means are contained e.g. in \cite{Triebel} and more
recents books by the same author.

 Strichartz 
\cite{Strichartz} found (see also \cite{Cho}) equivalent norms for 
$\dot H^{1,p}(\real^n)$ or, more generally, for corresponding spaces
with fractional smoothness and Lebesgue exponents in the range $p > n/(n+1)$. In this case the elements 
in 
spaces  $\dot H^{1,p}(\real^n)$ are locally integrable. Thus, $ \|f\|_{\dot H^{1,p}}\sim \| D_{2,1}f\|_{p}
,$ where $p >n/{n+1}$, and
$$
D_{2,1}(f) (x) =\left (  \int_0^\infty \left[ \int_{B(0,1)}|\Delta^2_{ry}f(x)\, dy|\right ]^2 r^{-3}\, dr\right )^{1/2}.
$$ 
Above $\Delta^2_{t} f(x)=f(x+2t)-2f(x+t)+f(x)$, whence Strichartz's characterization is pointwise, but it
employs integrated second differences.

Miyachi \cite{Mi3} characterized the Hardy-Sobolev spaces $\dot H^{1,p}$ in terms of
maximal functions related to  mean oscillation of the function in  cubes, thus
obtaining a
a counterpart of previous results of Calderon and of  the general theory of DeVore and Sharpley \cite{DeVoreSharpley}.
 More recently there has been considerable interest
in Hardy-Sobolev spaces $H^{1,p}$ and their variants on $\real^n$, or on subdomains. Chang, Dafni, and  Stein
\cite{ChangDafniStein} (see also \cite{ChangKrantzStein}) consider Hardy-Sobolev spaces in
connection with
estimates for elliptic operators, whereas  Aucher, Emmanuel, and Tchamitchian 
\cite{AucherEmmanuelTchamitchian} study these
spaces with applications to  square roots of elliptic operators.
Also the papers of Gatto, Segovia, and Jimenez \cite{GaSeJi}, Janson \cite{Janson} and  Orobitg \cite{Orobitg}
are related to the theme of the present paper.

Our main result shows that there is a  surprisingly simple {\it strictly pointwise}
characterization of the homogeneous (Hardy-)Sobolev space simply {\it in terms of  first differences:}
\begin{theorem}\label{th:main} Let $n\geq 1$ and ${n\over n+1} <p\leq 1.$ Then a distribution $f\in S'(\real^n)$
belongs to $\dot H^{1,p}(\real^n )$ if and only if $f$ is locally integrable and
there is a function $g\in L^p$ such that 
\begin{equation}\label{eq:1.10}
|f(x)-f(y)|\leq |x-y|(g(x)+g(y)),\quad x,y\in \real^n\setminus E
\end{equation}
where $E$ is a set of measure zero. Moreover, one has the equivalence of (quasi)norms
\begin{equation}\label{eq:1.12}
 \| f\|_{\dot H^{1,p}}\sim  \inf \| g\|_p,\nonumber
\end{equation}
where the infimum is taken over all admissible functions $g$ in {\rm \refeq{eq:1.10}}, and one identifies functions differing
only by a constant.
\end{theorem}
It was previously  known that the above characterization  holds true in the case $p>1.$ In \cite{Hajlasz1},
Hajlasz proposed to use \refeq{eq:1.10} as a definition of Sobolev-spaces on arbitrary metric
spaces. A considerable activity (see e.g.  \cite{Hajlasz1}, \cite{HajlaszKoskela},
\cite{Hajlasz2} and the references therein) has been devoted to the study of the corresponding (non-homogeneous) spaces
$M^{1,p}$, as they are customarily denoted. In the case $p>n/(n+1)$ we may define the non-homogeneus Sobolev
spaces $H^{1,p}(\real^n )$ by adding  to the definition of $\dot H^{1,p}(\real^n )$ the condition $f\in L^p(\real^n)$. When our result is combined with the previously known case $p>1$ (recall that
$L^p(\real^n)=H^p(\real^n)$ if $p>1$), we obtain the norm equivalence
\begin{equation}\label{eq:1.20}
 \| f\|_{H^{1,p}(\real^n)}\sim \| f\|_{M^{1,p}(\real^n)} \quad\mbox{for}\;\; p>{n\over n+1}, 
\end{equation}
which, incidentally,
solves the characterization problem of $M^{1,1}(\real^n)$ that has
been open after \cite{Hajlasz1}.
This fact also testifies for the naturality of the spaces of type $M^{1,p}$ in the case of Euclidean spaces: 
they yield the right spaces in view of harmonic analysis also in the case $p\leq 1.$

For the precise definitions of the Sobolev and Hardy spaces we refer to Section \ref{se:proof}, which also contains  auxiliary results and the proof of Theorem \ref{th:main}.
The remaining two sections   provide  examples
 of the  flexibility and strength
of Theorem \ref{th:main}. In the present paper we  aim 
 to concentrate   on key ideas, whence  we have not striven here for
most general results.

More specifically, Section \ref{se:domains} starts by treating the case where the derivative belongs to a local Hardy space. In addition, 
 pointwise characterizations analoguous to Theorem \ref{th:main} are given for the spaces $H^{1,p}(\Omega )$
defined on subdomains
$\Omega\subset\real^n$.  These results are obtained by slight modifications of the considerations of
Section \ref{se:proof}. Moreover, we scetch  a transparent proof of Jones-Miyachi's extension
result for uniform domains.

Section  \ref{se:hardyinequality} in turn applies our characterization to extend effortlessly the classical Hardy inequality also to case $p\leq 1$ in the framework of Hardy-Sobolev spaces. This is a novel range of exponents, since it
is well known that e.g. for the space
$W^{1,1}(\Omega )$ the Hardy inequality is not true even if the domain $\Omega$ is a ball.

\section{Definitions, auxiliary results and proof of Theorem \ref{th:main}}\label{se:proof}

We begin by shortly recalling  the relevant definitions and results from the theory of
real Hardy spaces. For the readers convenience we use as a principal reference the monograph
 \cite{Stein}.  
Fix a compactly supported function $\psi\in C_0^\infty (\real^n)$
with $\int \psi =1$ and with $\supp \psi\subset \{ x\; :\; |x|\leq 1\} .$
Assume that $p>0$ and consider $f\in \sdist (\real^n), $ i.e. $f$ is a tempered distribution. By the definition of Fefferman and Stein,  $f$ belongs to $H^p(\real^n)$
if and only if $\hlm_\psi f\in L^p(\real^n),$ where
$$
\hlm_\psi f(x):=\sup_{t>0}|f*\psi_t(x)|.
$$
Above $\psi_t(y)=t^{-n}\psi (y/t).$
The corresponding norm (quasi-norm in case $p<1$) is obtained by setting $\| f\|_{H^p}:=\| \hlm_\psi f\|_{L^p}.$ 
One may replace $\psi$ by any element in $\ssmooth (\real^n)$ and obtain an
equivalent norm. In what follows we shall denote by $\hlm f$ the standard Hardy-Littlewood maximal function of $f$.

We also recall a local version of the so-called grand maximal function. The basic significance of the grand maximal function comes from the fact that it allows one to use certain controlled families
of test functions in the maximal function, instead of a single one. For any radius $r>0$, natural
number $N\geq 1$, and any distribution $f\in \distr (B(x,r))$, we define
\beqla{eq:3.20}
\hlm_{r,N}f(x)=\sup_{\stackrel{\scriptstyle\varphi\in \smoothie{N}(B(x,r'))}{0<r'\leq r}} |\langle f,\varphi \rangle |,\nonumber
\eeq
and the class of functions $\varphi$ over which the supremum  is taken is 
$$
\smoothie{N}(B)=\{ \varphi\in C^\infty_0(B(x,r))\, :\, |\partial^\alpha \varphi|\leq r^{-n-|\alpha |}
\; \mbox{for}\;  |\alpha |\leq N \} .
$$
Observe that if $\varphi\in \smoothie{N}(B)$ then for any $t>0$ one has
$t^{-n}\varphi (t^{-1}\,\cdot )\in \smoothie{N}(tB)$. 
If the supremum in \refeq{eq:3.20} is taken over all positive radii, that is if $r=\infty$, we use the shorthand $\hlm_{N}f(x):=\hlm_{\infty,N}f(x).$
Assuming that $N$ is large enough (depending on $p,n$, see \cite[III 1.8, 5.9]{Stein}), a distribution  $f\in \sdist (\real^n) $ belongs to $H^p (\real^n)$ if and only if $ \hlm_{N}f(x)\in L^p(\real^n).$
In the range $p >n/(n+1)$ we may choose $N=1.$
For vector valued functions $f=(f_1,\ldots ,f_n)$ we set 
$\hlm_{r,N}f(x) =\max_{1\leq j\leq n}\hlm_{r,N}f_j(x).$

Recall that, in the case $p >1 $, the  homogeneous  Sobolev spaces $\dot W^{1,p}(\real^n)$ are defined by demanding that all the
first order distributional derivatives of $f$ lie in $L^p(\real^n).$

\begin{definition}\label{def:2.3}
Let $p>0$. We say that a tempered distribution $f$ on $\real^n$ belongs to 
the homogeneous (Hardy-)Sobolev space $\dot H^{1,p}(\real^n)$
if and only if  $D_jf\in H^p(\real^n)$ for each $j=1,\ldots ,n.$ Moreover
$$
\| f\|_{\dot H^{1,p}(\real^n)}:= \sum_{j=1}^n\| D_jf\|_{H^{p}(\real^n)},
$$
whence we obtain a (quasi)Banach space modulo constants. 
\end{definition}
\noindent The notation used above follows Triebel's \cite{Triebel} convention and it deliberately avoids confusion for $p=1,$ since obviously one has $\dot H^{1,p}(\real^n)=\dot W^{1,p}(\real^n)$
if $p>1,$ while for $p=1$ this breaks down.

In the present paper our main interest lies in the case where the elements in $\dot H^{1,p}$ are honest locally
integrable functions. This actually  happens for $p\geq n/(n+1)$. The following result is
well-known, but we obtain it as a corollary of our
proof of Theorem \ref{th:main}, see Remark \ref{re:2.0} below.

\begin{proposition}\label{pr:2.1} Assume that $n/(n+1)\leq p\leq 1$ and  let $p^*:=\ds {pn\over n-p}$ be the Sobolev
conjugate exponent, so that $p^*\geq 1.$ Then
 $\dot H^{1,p}(\real^n)\subset L^{p^*}_{loc}(\real^n).$ 
Especially, the elements of $\dot H^{1,p}(\real^n)$ are locally integrable.
\end{proposition}

We next recall the definition of the spaces $M^{1,p}(\real^n)$ and their homogeneous
counterparts $\dot M^{1,p}(\real^n)$.
\begin{definition}\label{def:2.5}
 Let $p>0$ and let $\Omega\subset\real^n$ be a subdomain. A measurable function $u$ belongs
to the homogeneous Sobolev space $\dot M^{1,p}(\Omega )$ 
(defined in the sense of Hajlasz) if there is a function $g\in L^p(\Omega )$ and a 
a set $E\subset\Omega $ of measure zero such that
for all $x,y\in \Omega \setminus E$ we have the estimate
\beqla{eq:2.40}
|u(x)-u(y)|\leq |x-y|(g(x)+g(y)).
\eeq
The corresponding quasi-norm is obtained by setting 
\begin{equation}\label{eq:2.50}
 \| f\|_{\dot M^{1,p}(\Omega )}:=  \inf \| g\|_{L^p(\Omega )},\nonumber
\end{equation}
where the infimum is taken over all admissible functions $g$ in \refeq{eq:2.40}. The  non-homo\-geneous
space $M^{1,p}(\Omega )$ is obtained by requiring, in addition, that $f\in L^p(\Omega )$, and the
norm for this space is defined by $ \| f\|_{M^{1,p}(\Omega )}:= \| f\|_{\dot M^{1,p}(\Omega )} + \| f\|_{L^p(\Omega )}.$
\end{definition}
Usually, we do not specifically mention the exceptional set $E$ since one may naturally allow the function $g$ to have value
$\infty$ in a set of measure zero. If \refeq{eq:2.40} holds for a certain measurable $g$
we  say that $g\in D(u).$
 It is well-known that for $p>1$ one has $M^{1,p}(\real^n)=W^{1,p}(\real^n)$. Actually, in this case 
(see \cite[Thm 2.2 and  formula (2.5)]{Hajlasz2}) there is a constant $c>0$ so that
\beqla{eq:2.55}
g:=c\hlm(|Df|)\in D(u).
\eeq
 When $p<1$ the spaces under consideration  are, of course,
quasi-Banach spaces (modulo constant functions), and for $p\geq 1$ Banach spaces. 
For these and other basic facts on
 $M^{1,p}(\real^n)$ we refer to \cite{Hajlasz2}. Especially, Hajlasz established an  important
extension \cite[Thm. 8.7]{Hajlasz2} of the Sobolev embedding  theorem to the quasi-Banach case, which holds true for metric spaces satisfying a lower bound
for the growth of the measure of balls.  In the case of $\real^n$ a special case of Hajlasz's 
theorem states the following.  
\begin{proposition}\label{pr:2.4}
Let $B\subset\real^n$ be a ball with radius $r$.
Then,  for any  $u\in \dot M^{1,p}(2B)$ with $p>0$, and $g\in D(u)$ one has the estimate
$$
\inf_{c\in\real}\left( \mint_B |u-c|^{p^*}\, dm\right)^{1/p^*} \leq Cr\left(\mint_{2B}g^{p}\, dm\right)^{1/p}.
$$
Here  $p^*={np\over n-p}$ and $C$ depends only on $p$ and $n$.
\end{proposition}
\noindent Above one may replace the domain of integration  $2B$ by any ball $\lambda B$, where $\lambda \in (1,2].$
If $p^*\geq 1,$ or equivalently $p\geq n/(n+1),$ we see that an element $f\in M^{1,p}(\real^n )$ is
locally integrable, whence it defines a distribution. A standard argument which uses 
Proposition \ref{pr:2.4} to compare
mean values of $f$ in balls $B_j$ and $B_{j+1}$ with $B_j=B(0,2^j)\subset\real^n$ shows that
$u(x)(1+|x|)^{-m}\in L^1(\real^n)$ for large enough $m.$ Hence $u$ (and, consequently, its derivatives) lie in 
$\sdist (\real^n).$

We will apply the previous Proposition in our proof of 
Theorem \ref{th:main}. The other essential ingredient is
Theorem \ref{pr:3.1} below, which provides the appropriate generalization of \refeq{eq:2.55} to the
case $p\leq 1.$  In the proof of Theorem \ref{pr:3.1} we need the following 
lemma, which is certainly well-known, but we were not able to find a suitable reference.
\begin{lemma}\label{le:smooth}
Let $Q\subset\real^n$ be a  cube and $\varphi\in C_0^\infty (Q).$ Then there
are elements $\psi_k\in C_0^\infty (Q)$, $k=1,\ldots ,n$ such that
$$
\varphi=\sum_{k=1}^nD_k\psi_k
$$
if and only if the condition $\int_{\real^n}\varphi \, dx =0$ is satisfied.
\end{lemma}
\begin{proof}
The stated condition is trivially necessary. In order to prove the sufficiency 
we apply induction on $n$. The case $n=1$ is evident, so suppose that
the result holds true for a fixed $n\geq 1$ and let $Q=Q'\times I\subset\real^{n+1}$, where $Q'\subset\real^n$ is a cube and
$I\subset \real$ is an interval. For $x\in\real^ {n+1}$
we 
 write $x=(x',x_{n+1})$ where $x'=(x_1,\ldots ,x_n)\in \real^n.$ 
Assume  that $\varphi\in C^\infty_0(Q)$
has zero mean. Define
$$ 
h(x'):=\int_{-\infty}^{\infty}\varphi (x',u)\, du.
$$
Then $h\in  C^\infty_0(Q')$ and, moreover, $h$ has zero mean.
The induction hypothesis enables us to write
$$
h(x')=D_1h_1(x')+D_2h_2(x')+\ldots + D_nh_n(x')
$$
with $h_j\in C^\infty_0(Q')$ for $j=1,\ldots ,n.$ Finally,
pick $a\in C^\infty_0(I)$ with $\int_I a(u)\, du =1$ and observe that the desired
decomposition  is obtained by choosing
$$
\psi_{n+1}(x):=\int_{-\infty}^{x_{n+1}}(\varphi (x',u)-a(u)h(x'))\, du
$$
and $\psi_j(x)= a(x_{n+1})h_j(x')$ for $j=1,\ldots ,n.$
\end{proof}

If $x,y\in\real^n$ we denote by $B_{x,y}$ the ball with the segment between $x$ and $y$ as 
a diameter. Observe that the only assumption on $f$ below is  the local integrability. The following
result is of independent interest.
\begin{theorem}\label{pr:3.1} For any $N\geq 1$ there exists a constant $c=c(N,n)$ with the following
property:
If $B\in\real^n$ is a ball and  $f\in L^1(2B)$,  then there is a set $E\subset B$ of measure zero such that for 
every $x,y\in B\setminus E$ it holds that
$$
|f(x)-f(y)|\leq c|x-y|(\hlm_{|x-y|,N}Df(x)+\hlm_{|x-y|,N}Df(y)).
$$
\end{theorem}
\begin{proof} 
Fix $N\geq 1.$ By rotational symmetry we may assume that $x-y=re_1$, where $e_1$
is the first unit coordinate vector and $r=|x-y|>0$.
 Let $\varphi\in C_0^\infty (B(0,1))$ be a fixed test function with $\int_{\real^n}\varphi\, dx =1.$
Choose  $k_0\geq 1$ so that $2^{k_0-1}\geq \sqrt {n}.$

Let us denote  
$$
A_k=\int_{\real^n} f(z+x)2^{n(k+k_0)}r^{-n}\varphi ({2^{k_0+k}r^{-1}}z)\,dz \quad \mbox{for} \; k\geq 0.
$$
 and analogously
$$
B_k=\int_{\real^n} f(z+y)2^{n(k+k_0)}r^{-n}\varphi ({2^{k_0+k}r^{-1}}z)\, dz\quad \mbox{for} \; k\geq 0.
$$
We apply Lemma \ref{le:smooth} to the function
$\varphi (\cdot )-2^n\varphi (2\, \cdot )$ and write
$$
\varphi (z)-2^n\varphi (2z)=\sum_{j=1}^nD_k\psi_j (z)\quad\mbox{for}\:\: z\in\real^n
$$
with $\psi_j\in \smooth ([-1,1]^n)$ for $1\leq j\leq n.$
An integration by parts yields that
\beqla{eq:3.27}
|A_k-A_{k+1}|&=&\big|\sum_{j=1}^n\int_{\real^n} f(z+x)
2^{n(k+k_0)}r^{-n}D_j\psi_j ({2^{k_0+k}r^{-1}}z)\,dz\big|\nonumber\\
&=& 2^{-k_0-k}r\big|\sum_{j=1}^n\langle D_jf(\cdot+x)2^{n(k+k_0)}r^{-n}\psi_k({2^{k_0+k}r^{-1}}\cdot )
\rangle\big|\nonumber\\
&\leq& C(n,N)r2^{-k}\hlm_{|x-y|,N}Df(x)\quad \mbox{for} \; k\geq 0.
\eeq
Similarly
\beqla{eq:3.28}
|B_k-B_{k+1}|
\leq C(n,N)r2^{-k}\hlm_{|x-y|,N}Df(y) \quad \mbox{for} \; k\geq 0.
\eeq

It remains to estimate the difference $|A_0-B_0|.$ Denote
 $$
\widetilde \varphi (z) = 2^{nk_0}[\varphi (2^{k_0}z)
-\varphi (2^{k_0}(z+e_1))].
$$ Then $\int_{\real^n}\widetilde\varphi\, dx=0$. By the choice of
$k_0$ we may apply Lemma \ref{le:smooth} to $\widetilde \varphi$ to obtain
functions $\widetilde\psi_k\in\smooth (B(0,1)\cup B(e_1,1))$ such that
$\widetilde\varphi=\sum_{k=1}^nD_k\widetilde\psi_k.$ With the help of suitable
cut-off functions we  may for each $k$ decompose  $k$ $\widetilde\psi_k= \widetilde\psi_{k,1}+\widetilde\psi_{k,2},$
where supp($\widetilde\psi_{k,1})\subset B(0,1)$ and supp($\widetilde\psi_{k,2})\subset B(e_1,1)$.
Then by translating, scaling and integrating by parts
we obtain as before that
\beqla{eq:3.29}
|A_0-B_0|\leq Cr(\hlm_{|x-y|,N}Df(x)+\hlm_{|x-y|,N}Df(y)).
\eeq

Finally, observe that if both $x$  and $y$ are   Lebesgue points of $f$ we have $f(x)=\lim_{k\to\infty}A_k$ and $f(y)=\lim_{k\to\infty}B_k$. The estimates \refeq{eq:3.27}--\refeq{eq:3.29} thus yield that 
\beqla{eq:3.30}
|f(x)-f(y)|&\leq& |A_0-B_0|+\sum_{k=0}^\infty (|A_k-A_{k+1}|+|B_k-B_{k+1}|\nonumber\\
&\leq& C'r(\hlm_{|x-y|,N}Df(y)+\hlm_{|x-y|,N}Df(y)),\nonumber
\eeq
and this finishes the proof of the Theorem.
\end{proof}

We are now ready for the proof of our main theorem.

\begin{proof}[Proof of Theorem \ref{th:main}.]
Assume first that $n/(n+1)<p\leq 1$ and $f$ is a Schwartz distribution on $\real^n$ such that $D_jf\in H^p(\real^n)$ for all $1\leq j\leq n.$
By Proposition \ref{pr:2.1} we have that $f\in L^1_{loc}(\real^n)$
(see also Remark \ref{re:2.0} below). Hence Theorem \ref{pr:3.1}
applies and we deduce that outside an exeptional set of measure zero
$f$ satisfies the inequality \refeq{eq:2.40} with the choice
$g=\hlm_{1}Df$. By the assumtions 
we have  $g\in L^p(\real^n)$. Moreover,
there is the estimate $\| g\|_{L^p}\leq C\sum_{j=1}^n\| D_jf\|_{L^p}.$
We have shown that $f\in \dot M^{1,p}(\real^n)$ with the 
correct bound for the (quasi)norm.

For the converse, assume next that $f\in \dot M^{1,p}(\real^n)$ and $p> n/(n+1).$ Fix $j\in\{ 1,\ldots ,n\} .$ We are to show that
$D_jf\in H^p(\real^n).$ For that end, we denote by $B(x,r)$ the ball of
radius $r$ and center $x$. Recall also that the support of $\psi$ is contained in the open unit ball,
and denote $C_0:=\|D\psi\|_\infty$. Since obviously $f_{|B(x,2r)}\in \dot M^{1,{n/(n+1)}}(B(x,2r))$, we may apply Proposition
\ref{pr:2.4} with $p=n/(n+1)$ and obtain that
 \beqla{eq:2.60}
\mint_{B(x,r)} |f-f_{B(x,r)}|\, dm \leq C_0r\left(\mint_{B(x,2r)}g^{n/(n+1)}\, dm\right)^{(n+1)/n},
\eeq
where $f_{B(x,r)}$ stands for the mean value of $f$ in the ball $B(x,2r)$. Apply the
above inequality to compute
\beqno
\hlm_\psi(D_ju)(x)&=& \sup_{t>0}|\langle D_jf,t^{-n}\psi ((x-\cdot)/t)\rangle |\\
&=& \sup_{t>0}|\langle f,t^{-n-1}(D_j\psi) ((x-\cdot)/t)\rangle |\\
&=& \sup_{t>0}|t^{-n-1}\int_{B(x,t)} (f(y)-f_{B(x,t})(D_j\psi) ((x-y)/t)\, dy|\\
&\leq& C_0\sup_{t>0}t^{-1}\mint_{B(x,t)} |f(y)-f_{B(x,t}|\, dy\\
&\leq& 2CC_0\sup_{t>0}\left( \mint_{B(x,2t)}g^{n/(n+1)}\right)^{(n+1)/n}\\
&\leq & C' (\hlm(g^{n/(n+1)})(x))^{(n+1)/n}.
\eeqno
By the Hardy-Littlewood theorem and the assumption
$g\in L^p(\real^n)$ we deduce that $\hlm (g^{n/(n+1)})\in L^q(\real^n)$ where $q=(n+1)p/n$. This shows that
$\hlm_\psi(D_iu)\in L^p(\real^n)$, as was to be shown.
\end{proof}

\begin{remark}\label{re:2.0}
Actually we may easily bypass the use of Proposition \ref{pr:2.1}
in the above proof. Namely, consider a Schwartz distribution $f$ such
that such that $D_jf\in H^p(\real^n)$ for all $1\leq j\leq n.$
Choose smooth convolution approximations $f_k$ of $f$ so that
$\|D_jf_k-D_jf\|_{H^p(\real^n)}\leq 2^{-k}$ for each $k\geq 1$ and 
$j=1,\ldots , n.$ By applying the above proof on each difference
$f_{k+1}-f_k$ and by writing $f=f_1+\sum_{k=1}^\infty (f_{k+1}-f_k)$
we obtain the desired result for $f$. In turn, Proposition \ref{pr:2.1}
is now obtained as a consequence of Theorem \ref{th:main} and 
Proposition \ref{pr:2.4}.

\end{remark}

Note that according to Proposition \ref{pr:2.1} we could as well assume a priori that $f$ is a tempered distribution
in the following definition.
\begin{definition}\label{def:2.7}
Let $p<n/(n+1)$. We say that a locally integrable $f$ on $\real^n$ belongs to the (non-homogeneous) Hardy-Sobolev space $H^{1,p}(\real^n)$
if   $D_jf\in H^p(\real^n)$ for each $j=1,\ldots ,n$ and $f\in L^p(\real^n).$ Moreover, we set
$$
\| f\|_{H^{1,p}(\real^n)}:= \sum_{j=1}^n\| D_jf\|_{H^{p}(\real^n)}+ \| f\|_{L^p(\real^n)}.
$$
\end{definition}

The following immediate corollary of Theorem \ref{th:main} verifies \refeq{eq:1.20}.
\begin{corollary}\label{cor:2.3}
Let $p>n/(n+1)$. Then $M^{1,p}(\real^n)=H^{1,p}(\real^n)$ with equivalent norms. 
\end{corollary}

\bigskip

It is well-known that bounded sets in the Hardy-spaces are weakly compact (see \cite[p. 127]{Stein}) in the
sense that any bounded sequence contains a subsequence that converges in the sense of distributions
to an element in the same space.
We thus obtain  the following compactness result for the spaces $M^{1,p}(\real^n)$.
\begin{corollary}\label{cor:2.1} Let $p>n/(n+1).$ Bounded sets in the  space $M^{1,p}(\real^n)$  are weakly
compact in the following sense: if $(f_k)$ is a norm-bounded sequence in $M^{1,p}(\real^n)$,
then there is an element $f\in M^{1,p}(\real^n)$ and a subsequence $f_{k_\ell}$ so that
$$
f_{k_\ell}\to f\quad \mbox{in}\;  L^1_{loc}\quad \mbox{as}\;\; \ell\to\infty .
$$
\end{corollary}
\begin{proof} Consider a norm-bounded sequence $(f_k)$ in $M^{1,p}(\real^n)$. By the above mentioned weak compactness of the Hardy spaces
we may pass to a subsequence and assume that for each $j$ it holds that $D_jf_k\to h_j$
in the sense of distributions. Obviously $h_j=D_jg$ for some distribution $g.$ Then $g\in H^{1,p}(\real^n)$ and the rest
follows by an application of Corollary \ref{cor:2.3} via the  known
compactness \cite{Triebel} of the embedding (restriction map) $H^{1,p}(\real^n)
\subset L^1(B)$ for any ball $B\subset\real^n.$
\end{proof}

Regarding applications to $H^{1,p}$ spaces we point out that the truncation stability
(or, more generally, Lipschitz stability) of these spaces,  proven by Janson \cite{Janson},
is an immediate consequence of our Theorem \ref{th:main}.

\begin{remark}\label{rem:2.2} One may check that
$H^{1,p}(\real^n)=\dot F^{1}_{p,2}(\real^n)\cap L^p(\real^n)$, where $\dot F^{1}_{p,2}(\real^n)$ is the homogeneous
Triebel space, see \cite[Section 5]{Triebel}. 
\end{remark}

\begin{remark}\label{rem:2.5} 
In  the proof of Theorem \ref{th:main} it would be possible to partially apply
the existing results on 
Hardy-Sobolev spaces \cite{GaSeJi},  \cite{Janson},  \cite{Mi3}, \cite{Orobitg}.
However, our proof is direct and simple and it employs just the maximal function
definition of the Hardy spaces.
It is perhaps of interest to also note
that one may also prove first half of Theorem \ref{th:main}, i.e. the inclusion $\dot H^{1,p}(\real^n)\subset \dot M^{1,p}(\real^n)$
by applying the atomic decomposition of the Hardy spaces
(see \cite[III.2]{Stein}). This is done with the aid of the representation
$$
f=\sum_{j=1}^nI_1R_jD_jf +const,
$$
where the $R_j$:s  are the Riesz transforms and the Riesz potential $I_1$  corresponds to the Fourier-multiplier
$\widehat f\mapsto |\xi |^{-1} \widehat f.$ Since the Riesz transforms
are bounded on the Hardy spaces, and we are dealing with the case $p\leq 1$, it turns out that it is enough to show that for
each $H^p$-atom $a$ we have $I_1a\in \dot M^{1,p}(\real^n).$
Moreover, by scaling and translation invariance one may assume that 
the atom $a$ is related to the unit ball $B(0,1)$, whence one
may prove by hand that $I_1a\in \dot M^{1,p}(\real^n).$
On the other hand, Theorem \ref{pr:3.1} is a more convenient tool
since it bybasses the atomic theory and, more importantly, it
is a local result which applies directly to spaces 
defined on  subdomains of $\real^n,$ see Section \ref{se:domains}
below. 
\end{remark}

\begin{remark}\label{rem:2.3} It is natural to also consider
exponents $p\leq n/(n+1)$ and ask for the right analogue for our characterization
of the Hardy-Sobolev spaces. Moreover, we do not know if Theorem
\ref{th:main} holds true as such in the case $p=n/(n+1).$
\end{remark}

\section{Local spaces and spaces on subdomains of $\real^n$}\label{se:domains}

The elements of $H^p(\real^n)$-spaces satisfy moment conditions, e.g. any function from $H^1(\real^n)$ has mean zero.
In order to relax this condition and to obtain localizable  spaces, one defines (\cite{Goldberg}, see
also \cite[Ch. 3, Sec. 5.17]{Stein}) the {\sl local} Hardy space $h^p(\real^n)$  by restricting, in the definition
of the maximal operator $\hlm_\psi$, the range of $t$ to the interval $(0,\rho ]$, where $\rho >0$ is a fixed positive number.
Define
$$
\hlm_{\psi,\rho} f (x):= \sup_{0<t\leq \rho}| f* \psi_t (x)|,
$$
whence a tempered distribution is said to belong to $h^p(\real^n)$ ($p>0$) if and only if $\hlm_{\psi,\rho} f\in L^p(\real^n).$ One
sets
 $\| f\|_{h^p}:=\| \hlm_{\psi,\rho} f\|_{L^p}.$
The above definition does not depend on $\rho >0$ or on $\psi .$ In particular,
different values of $\rho $ lead to equivalent (quasi)norms.
Finally, in terms of grand maximal functions   
$$
f\in h^p(\real^n)\quad\mbox{if and only if } \quad \hlm_{\rho,N}f\in L^p(\real^n).
$$
as soon as $N\geq N_0(p,n)$.

The local Hardy spaces   form a more flexible class than the Hardy spaces. There are no more global
moment conditions,
the elements in $h^p$ can be localized and they are (at least locally) invariant under diffeomorphisms.

 {\it The local Hardy-Sobolev space} is  obtained in a natural way:
\begin{definition}\label{def:3.3}
Let $p>0$. We say that a tempered distribution $f$ on $\real^n$ belongs to 
the homogeneous Hardy-Sobolev space $\dot h^{1,p}(\real^n)$
if and only if  $D_jf\in h^p(\real^n)$ for each $j=1,\ldots ,n.$ Moreover
$$
\| f\|_{\dot h^{1,p}(\real^n)}:= \sum_{j=1}^n\| D_jf\|_{h^{p}(\real^n)}.
$$
In case $p<1$ the above space is  a quasi-Banach space modulo constants. 
Similarly, in case $p\geq n/(n+1)$ we say that $f\in  h^{1,p}(\real^n)$ if $D_jf\in h^p(\real^n)$ for each $j=1,\ldots ,n,$
and $f\in L^p(\real^n).$
\end{definition}
\noindent The last sentence above makes sense since locally the elements of $h^p$ coincide with elements in $H^p$ whence Proposition \ref{pr:2.1} 
easily implies  $\dot h^{1,p}(\real^n)\subset L^1_{loc}(\real^n).$

We now verify a counterpart of Theorem \ref{th:main} for the local spaces. Below one may replace the condition
$|x-y|\leq 1$ by $|x-y|\leq c$, where $c>0$ is an arbitrary constant.
\begin{theorem}\label{th:3.1} Let $n\geq 1$ and ${n\over n+1} <p\leq 1.$ Then a distribution $f\in S'(\real^n)$
belongs to $\dot h^{1,p}(\real^n )$ if and only if $f$ is locally integrable and
there is a function $g\in L^p$ such that 
\begin{equation}\label{eq:3.200}
|f(x)-f(y)|\leq |x-y|(g(x)+g(y)),\quad\mbox{for}\; \; |x-y|\leq 1\;\; \mbox{and}\; x,y\in \real^n\setminus E
\end{equation}
where $E$ is a set of measure zero. Moreover, one has the equivalence of (quasi)norms
\begin{equation}\label{eq:3.201}
 \| f\|_{\dot h^{1,p}}\sim  \inf \| g\|_p,\nonumber
\end{equation}
where the infimum is taken over all admissible functions $g$ in {\rm\refeq{eq:3.200}}, and one identifies functions differing
only by a constant.
\end{theorem}
\begin{proof}
Let us first assume that condition \refeq{eq:3.200} holds true. If we replace $\sup_{t>0}$
by $\sup_{t\leq 1/4}$ in the second part of the proof of Theorem \ref{th:main} we may again
apply Proposition \ref{pr:2.4} in a similar manner to deduce that
$
\hlm_{\psi,1/4} (D_kf)(x)\in L^p(\real^n).
$
It follows that $f\in \dot h^{1,p}(\real^n )$. 

In order to prove the converse, observe that according to Theorem \ref{pr:3.1} condition
\refeq{eq:3.200} is satisfied with the choice
$
g:=c\max_{1\leq j\leq n}\hlm_{N,1} D_jf .
$
The assumption $D_jf\in h^p$ now implies that $g\in L^p$ as soon as
$N$ is large enough. 
Here, as also in the first part of the proof, the corresponding  quantitative
statement is obvious.
\end{proof}

We now turn to the study of Hardy-Sobolev spaces on subdomains of the Euclidean space. A
natural and simple definition of Hardy spaces 
that works for all subdomains $\Omega\subset\real^n$ was given by Miyachi \cite{Mi1}.
We recall a definition that is directly equivalent with Miyachi's definition. Thus, let $\Omega\subset\real^n$ be a domain,  $p\in (0,1)$ and $N>\max (0,n({1\over p}-1))$.  If 
 $f\in \distr (\Omega )$ one  defines
\begin{equation}\label{eq:3.210}
f\in H^p(\Omega ) \quad \mbox{if and only if} \quad
\| f\|_{H^p(\Omega )}=:\| \hlm_{d(x,\boundary\Omega )/2,N}f\|_{p}<\infty
\end{equation}
It follows from Miyachi's results that the  above is equivalent with the original definition \cite{Mi1}. Similarly, we see immediately that   one may replace above $\hlm_{d(x,\boundary\Omega )/2}f$ by
$\hlm_{ad(x,\boundary\omega )}f$ for any $a\in (0,1)$, or by a the more standard maximal
function
$\sup_{0< t\leq d(x,\boundary\Omega )}|f*\psi_t|,$ where the fixed test function
$\psi\in C_0^\infty(B(0,1))$ has nonzero mean. All the (quasi)norms so obtained are mutually
equivalent.
Let us also note that prior to \cite{Mi1}, Jonsson, Sj\"ogren, and Wallin \cite{JoSjWa} defined Hardy
spaces on fairly general subsets  of $\real^n$ in terms of suitable
atoms.

Recall that the spaces $\dot M^{1,p}(\Omega )$ and $ M^{1,p}(\Omega )$
were defined already in the previous section. 
The spaces $\dot M^{1,p}_{ball}(\Omega )$ and $M^{1,p}_{ball}(\Omega )$
are defined exactly in the same manner but for one difference: in the definition of these
spaces the condition \refeq{eq:2.40} is assumed to hold only for points
$x,y\in\Omega\setminus F$ that satisfy the condition 
$$
|x-y|\leq {1\over 4}\min (d(x,\boundary\Omega ),
d(y,\boundary\Omega )).
$$
After the proof of Theorem \ref{th:3.4} below it is clear that one may replace in the last condition above the constant
$1/4 $ by any constant strictly less than one.

The Hardy-Sobolev spaces
on $\Omega$ are defined in the obvious manner.
\begin{definition}\label{def:3.2} Let $\Omega\subset\real^n$ be a domain and
$p>0.$ A measurable function $u:\Omega\to\real$ belongs
to the homogeneous Hardy-Sobolev space $\dot H^{1,p}(\Omega )$ 
 if $D_jf\in H^p(\Omega )$ for all $j=1,\ldots ,n.$
The related seminorm is obtained by setting 
$
 \| f\|_{\dot H^{1,p}(\Omega )}:=  \sum_{j=1}^n \| D_jf\|_{H^p(\Omega )} $.
In case $p\geq n/(n+1)$ we say that $f\in  H^{1,p}(\Omega )$ if
additionally  $f\in L^p(\Omega )$, and the corresponding (quasi)norm
is defined in the obvious manner. 
\end{definition}

\begin{theorem}\label{th:3.3} Let $n\geq 1$ and ${n\over n+1} <p\leq 1$, and assume that
$\Omega\subset\real^n$ is a domain. Then 
$$
H^{1,p}(\Omega )= M^{1,p}_{ball}(\Omega )\quad \mbox{and}\quad 
\dot H^{1,p}(\Omega )= \dot M^{1,p}_{ball}(\Omega ),
$$
with equivalence of the (quasi)norms.
\end{theorem}
\begin{proof}
It is enough to prove the latter equality.
Assume first that $f\in H^{1,p}(\Omega).$ Let  $|x-y|\leq {1\over 4}\min (d(x,\boundary\Omega ),
d(y,\boundary\Omega ))$. Since the first derivatives of $f$ are
locally in $h^{1,p},$ it easily follows that $f\in L^1_{loc}$. We may apply Theorem \ref{pr:3.1}
in the ball $B({1\over 2}(x+y),{1\over 3}d(x,\partial\Omega))$
in order to obtain the inequality \refeq{eq:2.40} with the choice
$$
g(x)= \hlm_{{1\over 2}d(x,\partial\Omega), N}Df(x).
$$
By definition it holds that $g\in L^p(\Omega ).$

The converse follows exactly as the proof of the second part of Theorem \ref{th:main}.
One just uses the observation that 
$
\sup_{\varphi\in \smoothie{N}(B(x,r))}\|D_k\varphi\|_1\leq c(n,N)r^{-1}
$ for each $k=1,\ldots, n.$
\end{proof}

For simplicity,  from now on we state our results only for non-homogeneus spaces. The
reader will have no difficulty in formulating the corresponding results for
the homogeneous spaces.
After our previous results it is of interest to find conditions on $\Omega$ that would  quarantee
that the obvious inclusion $ M^{1,p}(\Omega )\subset M^{1,p}_{ball}(\Omega )$ becomes
equality. Let us recall for that end the definition of {\it uniform domains}. One says that a domain
$\Omega\subset\real^n$ is uniform if there is a constant $c>0$ such that 
for all $x,y\in \Omega$ there is a rectifiable path $\gamma: [0,T]\to \Omega$, parametrized
by arclength,  with $\gamma(0)=x$, $\gamma(T)=y$, and such that $T\leq c|x-y|$ together with
$$
B(\gamma (t), {1\over c}\min (t,T-t))\subset\Omega \quad \mbox{for}\;\;
t\in (0,T).
$$
\begin{theorem}\label{th:3.4} Let $n\geq 1$ and ${n\over n+1} <p\leq 1$, and assume that
$\Omega\subset\real^n$ is a uniform domain. Then $M^{1,p}_{ball}(\Omega )=M^{1,p}(\Omega )$
and, especially
$$
H^{1,p}(\Omega )= M^{1,p}(\Omega )
$$
with equivalent (quasi)norms.
\end{theorem}
\begin{proof} 
Choose arbitrary $x,y\in \Omega .$ From the uniformity condition we deduce easily that
there is a chain of balls $B_k$ resembling a cigar  that joins the points $x$ and $y$.
To be more specific, there are balls $B_k:=B(z_k,r_k)$ with $k\in\integer$ and $z_k\in \Omega$ such that for each $k$ one has 
$6B_k\subset \Omega $,  $r_k\geq {1\over c'}\min (d(z_k,x),d(z_k,y) ),$
$B_k\cap B_{k+1}\not=\emptyset $,  and $r_k/2\leq r_{k+1}\leq 2r_k.$
In addition, $\lim_{k\to\infty} 
d(x,B_k)=0=\lim_{k\to-\infty} d(y,B_k).$
Finally, we may assume that $\sum_{k\in\integer}r_k\leq c'|x-y|$

Let $f\in M^{1,p}_{ball}(\Omega )$. We  compare in a standard manner the consequtive
mean values   $f_{B_k}$ by applying the Poincare inequality of Proposition \ref{pr:2.4}
in the ball $3B_k$
 \beqla{eq:3.300}
|f_{B_k}-f_{B_{k+1}}| &\lesssim & \mint_{3B_k}|f-f_{3B_k}|\, dx
\lesssim  r_k\left(\mint_{6B_k}g^{n/(n+1)}\, dm\right)^{(n+1)/n}\nonumber\\
&\lesssim &  r_k((\hlm g^{n/(n+1)} (x))^{(n+1)/n}
+ (\hlm g^{n/(n+1)} (y))^{(n+1)/n}).\nonumber
\eeq
Above we assumed $g$ to be continued as zero outside $\Omega ,$
and the last estimate followed by observing that one of
the points $x,y$ is contained in $c'B_k.$

If $x,y$ are Lebesgue points of $f$ we have 
$|f(x)-f(y)|\leq \sum_{k\in\integer}|f_{B_k}-f_{B_{k+1}}|.$ By summing over $k$ in our
previous estimate it follows that 
$$
|f(x)-f(y)|\leq |x-y|(\widetilde g(x)+\widetilde g(y)),
$$
where $\widetilde g=c'C' ((\hlm g^{n/(n+1)})^{(n+1)/n})_{|\Omega}$.
 The conclusion
follows  from Theorem \ref{th:3.3} since by the standard Hardy-Littlewood maximal inequality
$\widetilde g\in L^p(\Omega ).$
\end{proof}

We next consider extensions from bounded domains to the whole  space $\real^n .$ We do not
aim for most general results here (see Remark \ref{rem:3.1} below), but show how our results
can be used to give a transparent new proof of the following result
of Miyachi \cite{Mi2}.
\begin{theorem}\label{th:3.5} (Miyachi) Let $n\geq 1$ and ${n\over n+1} <p\leq 1$, and assume that
$\Omega\subset\real^n$ is a bounded uniform domain. Then there is a bounded linear extension
operator from $H^{1,p}(\Omega )$ into $H^{1,p}(\real^n )$. 
A fortiori,
$$
H^{1,p}(\Omega ) = \{ f_{|\Omega } \; :\; f\in H^{1,p}(\real^n)\}.
$$
\end{theorem}
The above Theorem is a direct consequence of Theorems \ref{th:3.1} and \ref{th:3.4} as soon as
we verify  the following Proposition.
\begin{proposition}\label{pr:3.3} Let $n\geq 1$ and ${n\over n+1} <p\leq 1$, and assume that
$\Omega\subset\real^n$ is a bounded uniform domain. Then there is a bounded linear extension
operator from $M^{1,p}(\Omega )$ into $H^{1,p}(\real^n )$. 
A fortiori,
$$
H^{1,p}(\Omega ) = \{ f_{|\Omega } \; :\; f\in H^{1,p}(\real^n)\}.
$$
\end{proposition}
\begin{proof} \quad
Assume that $f\in M^{1,p}(\Omega)$ and $g\in L^p(\Omega )$ so that the condition \refeq{eq:2.40}
is satisfied. We assume that $g$ is continued as zero outside $\Omega$ and set $ g_1:=
(\hlm g^{\widetilde p})^{1/\widetilde p},$ where as before $\widetilde p$ is chosen
from the interval $(n/(n+1),p)$ so that $g_1\in L^p(\real^n)$. Assume that
$\overline{B}(x,r), \overline{B}'(x',r')\subset \Omega$ are two closed balls (or points).
 such that $2r\leq d(x,\boundary\Omega )$ and
$2r'\leq d(x',\boundary\Omega )$.
An easy modification of the proof 
 of Theorem \ref{th:3.4} yields the estimate:
\begin{equation}\label{eq:3.330}
|f_B-f_{B'}|\leq c(|x-x'|+r+r')(\inf_{z\in B}g_1(z)+\inf_{z'\in B'}g_1(z')).
\end{equation}
If, for example, $r=0$, we replace $f_B$ above by $f(x)$ and the corresponding infimum in the
right hand side by $g_1(x).$ This makes sense as long as $x$ is a Lebesgue point of $f.$ 

For $\varepsilon >0$ set $A_\varepsilon := \{ x\in\real^n\setminus\Omega \; :\; d(x,\Omega)\leq \varepsilon \}.$
We may choose $\varepsilon_0 >0$ and a Whitney type cover of
a neighbourhood of $\Omega$ by balls $\{ B_i\}_{i=1}^\infty$ so that for the set
 $A:=\cup_{i=1}^\infty B_i$ it holds that $A_{4\varepsilon_0}\subset A\subset A_{10\varepsilon_0}$. Moreover,
for each $i$ the balls $B_i=:B(y_i,\rho_i)$ satisfy 
$d(y_i,\Omega )\leq 4\rho_i\leq 2d(y_i,\Omega )$ and  the inflated balls $2B_i$ have bounded overlap.
Finally, the situation can be arranged so that we may also pick a partition of unity that consist of functions $h_i$
such that $\supp h_i\subset B_i ,$ and $0\leq h_i\leq 1$ together with $\| Dh_k\|_\infty\leq c/\rho_i.$ 
The desired partition of unity property may be expressed as follows:
$\chi_{A_{2\varepsilon_0}}\leq\sum_i h_i \leq\chi_{ A_{10\varepsilon_0}}.$

For each ball $B_i$ we may, by the uniformity of $\Omega$, select a ball $B'_i:=B'_i(y'_i,r'_i)\subset\Omega$
that satisfies uniformly $\diam (B_i)\sim \diam (B'_i)\sim  d(y'_i,y_i)\sim d(B'_1,\partial\Omega ).$
The required extension $F$ of $f$ to the whole space $\real^n$ is then simply
 defined by choosing $F(x)=f(x)$ for $x\in\Omega ,$ and 
$$
F(x)=\sum_{i=1}^\infty (\mint_{B_i}f\, dx)h_i(x)\quad \mbox{for  }\; x\in\real^n\setminus\overline{\Omega}. 
$$
Observe that by the uniformity of $\Omega$  we may ignore $\partial\Omega$ since 
it is of measure zero.

We claim  that condition \refeq{eq:1.10} is satisfied for a.e.  $x,y$ 
 with $g_2$ in place of $g$, where $ g_2:=
(\hlm g_1^{\widetilde p})^{1/\widetilde p}\in L^p(\real^n).$
The verification of this is a pretty routine using \refeq{eq:3.330}, and hence we just outline  
it. The details
will possess no difficulties for the reader.

Consider  arbitrary  $x,y\in A_{\varepsilon_0}\cup\Omega .$
We claim  that condition \refeq{eq:1.10} is satisfied for a.e. such points $x,y$. 
For $x,y\in\Omega $ this is evident, so assume that $x\in A_{2\varepsilon_0}\setminus\overline{\Omega},$
say $x\in B_i=B(y_i,r_i).$ The argument will be divided into different cases depending on
the relative location of the point $y.$ Observe that $F(x)$ is a convex combination of
averages of $f$ over   balls $B'_{i_1},\ldots , B'_{i_\ell}\subset\Omega $, say, whose distance to the boundary
$\partial\Omega$ is comparable to their size and to theír distance to $x.$ Moreover,
$\ell\leq C$ uniformly.

If $y\in\Omega ,$ it is thus enough to estimate differences of the form $|f(y)-f_{B'_{i_k}}|$.
We apply \refeq{eq:3.330} and obtain immediately an estimate of the desired type,
since in the present situation $g_2(x)$ dominates the quantity $\inf_{z\in B'_{i_k}}g_1(z)$,
and $|x-y|\gtrsim d(x,\cup_{k=1}^\ell f_{B'_{i_k}}).$

Let then also  $y\in A_{2\varepsilon_0}\setminus\overline{\Omega} .$ whence the value $F (y)$ is likewise computed
as a convex combination over averages in  certain balls of the cover. 
We consider first the case where $|x-y|\geq c'\max (d(x,\partial\Omega ),d(y,\partial\Omega )).$
The value $F (y)$ is likewise computed
as a convex combination over averages in  certain balls of the cover. 
If $B'_j=B(y'_j,r'_j)$ is one of these
balls, it is enough to estimate  the difference $|f_{B'_i}-f_{B'_j}|.$  This
is done essentially as in the previous paragraph again by applying \refeq{eq:3.330}.

The remaining possibility is that  $|x-y|\leq c'\max (d(x,\partial\Omega ),d(y,\partial\Omega )).$ 
If $c'$ is chosen appropriately, the properties of the Whitney type cover imply that distance of both $x$ and $y$ from the boundary $\partial\Omega$
are comparable to $r:=d(x,\partial \Omega)$ and the radii of the balls $B_k$ that contain one of the points $x,y$ are
 comparable to $r$, and the distance between these balls is less or comparable to $r$. Assume that $\rho_{k_1},\ldots ,\rho_{k_m}$ are the functions of the 
partition of unity that are nonzero at $x$ or at $y.$ Denote $a_u:=\mint_{B'_{k_u}}f_k.$
The estimate \refeq{eq:3.330} can be used to check that 
$$
B:=\max_{1\leq u,u'\leq m}|a_u-a_u'|\leq c''r(g_2(x)+g_2(y)).
$$
Hence, as the $h_i$:s form a partition of unity, we may estimate
$$
|F(x)-F(y)|\leq B\sum_{u=1}^m|h_{k_u}(x)-h_{k_u}(y)|.
$$
The desired upper bound for $ |F(x)-F(y)| $ now follows by combining the previous
estimates with the bound $|h_{k_u}(x)-h_{k_u}(y)|\leq |x-y|\|Dh_{k_u}\|\leq C|x-y|/r.$

We have shown that $F\in \dot M^{1,p}(A_{2\varepsilon_0}\cup\overline{\Omega})$. Finally, by
combining this with the fact that the function $F$ has compact support and is Lipschitz outside  
the set $A_{\varepsilon_0}\cup\overline{\Omega}$, we finally
obtain \refeq{eq:1.10} for a.e. $x,y$ in $\real^n$ 
e.g. with the function $c(g_2+F)$ substituted in place of $g$, where $c$ is a large enough
constant.  The estimate for $\| F\|_{L^p(\real^n)}$ is immediate and we obtain that
$F\in  M^{1,p}(\real^n)$ with suitable bounds for the norm.
\end{proof}

\begin{remark}\label{rem:3.1} The  motivation for presenting proof of Theorem \ref{th:3.5}
above was to demonstrate how the coincidence of the Sobolev-Hardy spaces
with the spaces $M^{1,p}$ can also be used to give simple and unified proofs of extension results, since
the above argument works unchanged also for $p>1$.
If one considers extension results just for the spaces  $M^{1,p}(\Omega )$ it is possible to
considerably weaken the conditions on $\Omega$ by replacing uniformity by
a so called measure density property of the domain $\Omega .$
 This result (in the case $p=1$) in contained in 
\cite{HaKoTu}.
\end{remark}

\section{Extension of the Hardy inequality}\label{se:hardyinequality}

For a subdomain $\Omega\subset\real $ we denote by $W^{1,p}_0(\Omega )$ the closure of $C_0^\infty (\Omega )$ in the space
$W^{1,p}(\Omega )$.
It is well-known that if $\partial \Omega$  is regular enough, then the classical Hardy inequality
\beqla{eq:4.10}
\int_{\Omega}\left({|u(x)|\over d(x,\partial\Omega)}\right)^p\,dx\leq  C(p,n,\Omega )
\int_{\Omega}|Du(x)|^p\,dx
\eeq
holds for $p\in (1, \infty ).$ In this case we say that $\Omega$ carries the $p$-Hardy inequality. We refer e.g.  to \cite{Lewis} for  conditions on $\Omega$
which ensure the validity of \refeq{eq:4.10}. On the other hand, it is well-known that the Hardy
inequality fails for the space $W^{1,1}_0$. This fact is closely connected to  
 $W^{1,1}(\real^n)$ and $M^{1,1}(\real^n)$ being different, as the 
following simple example shows.

\begin{example}\label{ex:4.1}
Assume that $h\in L^1(\real)$  satisfies supp$(h)\subset [-2,2]$, $h\geq 0$, and $h$ is even.
Define a function $u\in W^{1,1}(\real )$ by setting $u(x)=\int_0^{x}h(t)\, dt.$ Pick any
  $g\in D(u)\cap L^1(\real )$. By symmetry we may assume 
 that $g$ is even without changing the integrability properties of $g$. By
considering  the condition \refeq{eq:2.40} with respect to points $x$ and $-x$ we deduce that
\beqla{eq:4.0}
\big|{u(x)\over x}\big|\leq g(x)\quad \mbox{for a.e}\;\; x>0.
\eeq
Thus $u(x)/x$ is integrable in a neighbourhood of the origin, assuming that $g\in L^1.$ Especially this holds if $u\in M^{1,1}(\real ).$
In this case, \refeq{eq:4.0} can be viewed as  an analog of
the Hardy-inequality for the function $u.$ However, in general \refeq{eq:4.0}  fails for $u\in W^{1,1}(\real )$, since the choice
$h(x)\sim (|x|\log^2(1/|x|))^{-1}$ near the origin leads to $u(x)/x\sim (|x|\log(1/|x|))^{-1}$ for
small $x.$
\end{example}

We say that a bounded subdomain $\Omega\subset\real^n$ is Lipschitz if each point $x\in \partial\Omega$ has
a regular neighbourhood $U$ in the sense that,  after an isometric change of coordinates, we have
\beqla{eq:4.20}
U\cap\partial \Omega =\{ (y,h(t))\; :\; y\in B(0,\delta)\cap\real^{n-1}\}
\eeq
where $\delta >0$  and the Lipschitz function $h:B(0,\delta)\cap\real^{n-1}$ may depend on $x$, but the 
Lipschitz constants of the functions $h$ are uniformly bounded. This enables us to define a local
reflection $H$ on the neighbourhood $U\cap\Omega$. Thus, assuming that in the above coordinate system $\Omega$
lies locally below the graph of the function $h$, we define for
$x=(y,t)\in U\cap\Omega$ 
$$
H(x)= (y,2h(y)-t).
$$
The reflection map $H$ is a measure preserving bijection $H: U\cap\Omega\to H(U\cap\Omega )$. Morever, by shrinking the neighbourhood $U$
if needed, we may assume that  $\sim |H(x)-x| \leq c_1d(x,\partial\Omega)$,  since $h$ Lipschitz.

The following result shows that the right counterpart of the Hardy  inequality in the case $p\leq 1$ involves the Hardy-Sobolev spaces.
It should be noted that the proof is almost trivial thanks to Theorem \ref{th:main}.

\begin{theorem} Let $\Omega$ be a bounded Lipschitz subdomain of $\real^n$ and $p\in (n/(n+1),1]$. Then there is a finite constant
$C=C(p,n,\Omega)$ such that any $\varphi\in  C_0^\infty (\Omega )$ satisfies the inequality
\beqla{eq:4.40}
\int_{\Omega}\left({|u(x)|\over d(x,\partial\Omega)}\right)^p\,dx\leq  C\sum_{j=1}^n\|D_j f\|^p_{H^p(\real^n)}.
\eeq
\end{theorem}
\begin{proof}
Assume that $\sum_{j=1}^n\|D_j u\|^p_{H^p(\real^n)}=1$.
Choose a finite cover of of $\partial\Omega $ 
by regular neighbourhoods.
Let $U$ be one of these neighbourhoods and assume that its intersection with the
boundary $\partial\Omega$ has the representation \refeq{eq:4.20}. By Theorem \ref{th:main} we may choose
a function $g\in D(u)$ with $\| g\|_{L^p(\real^n)}\leq c.$ We note that $u$ vanishes outside $\Omega ,$ whence
we obtain for almost every $x\in U\cap\Omega$ that
$$
|u(x)|=  |u(x)-u(H(x))|\leq |x-H(x)|(g(x)+g(H(x)))\leq c_1 d(x,\partial\Omega )(g(x)+g(H(x)).
$$ 
Since $H$ is measure preserving we thus obtain that
$$
\int_{U\cap\Omega}({|u(x)|\over d(x,\partial\Omega)})^p\,dx\leq (2c_1)^p\int_{U\cap\Omega}g(x)^p\,dx \leq (2c_1)^pc^p.
$$
The claim follows   by summing over the chosen regular neighbourhoods, and applying the Poincare 
inequality to the uncovered part of the domain.
\end{proof}

\begin{remark}\label{rem:4.1}
Another, even easier  approach to the above result would be to use the obvious bi-Lipschitz
invariance of the space $M^{1,p}(\real^n)$ in order to reduce  the boundary
locally to a piece of a hyperplane, in which case the reflection argument is trivial.
\end{remark}

As our last observation we generalize the above observation to a more extensive class of
domains. 
For that end we recall certain definitions.
Let $U\subset\real^n$ be a bounded open set and $p\in (0,1].$  The   {\it Hardy $p$-capacity}
of a compact subset $E\subset U$, relative to $U$, is defined by setting
\beqla{eq:4.100} &&
\capa_{H^{1,p}}(E;U):=\inf \big\{ \sum_{k=1}^n\| D_k\varphi\|^p_{H^p(\real^n )}\; :\; \varphi\in C^\infty_0
(U),\; \varphi (x)\geq 1\; \mbox{for}\; x\in E\big\}.
\eeq
It is easily verified that one may replace the condition $\varphi\in C^\infty_0
(U)$ above by $\varphi\in W^\infty_0
(U)$, i.e. the class of compactly supported Lipschitz functions.

We need an analogue of \cite[Proposition 3]{Hajlasz3}. The proof is a slight modification of the original argument in case $p>1$.
\begin{lemma}\label{le:4.1} Let $q>n/(n+1)$. Assume that $B\subset\real^n$ is a ball with radius $r$,
 $K\subset B$ is compact with
$\capa_{H^{1,q}}(K;2B)\geq cr^{n-q}$. Then  for every $u\in C^\infty_0 (\real^n)$ satisfying $u_{|K}=0$, any $g\in D(u)$, and for any $x\in \overline{B}$
it holds that
$$
|u(x)| \leq rc(n,q)(M(g^q)(x))^{1/q}.
$$
\end{lemma}
\begin{proof}
For $k\geq 0$ denote $B_k=B(x,r2^{-k})$. The Poincare inequality, i.e. Proposition \ref{pr:2.4} yields  that
$|u_{2B}-u_{B_0}|\lesssim r(M(g^{q})(x))^{1/q}$ and 
$|u_{B_k}-u_{B_{k+1}}|\lesssim 2^{-k}r(M(g^{q})(x))^{1/q}$. Since $u$ is continuous we may sum up
these inequalities and obtain
\beqla{eq:4.1011}
|u(x)-u_{2B}|\lesssim r(M(g^{q})(x))^{1/q}.
\eeq
In order to estimate $u_{2B}$ we assume that $u_{2B}\not=0$. Define
a cutoff function $h$ by setting $h(x)=1$ on $B$, $h(y)=1-2{\rm dist}
(x,B)/r$ for $x\in (3/2)B\setminus B$, and $h=0$ for $x\not\in (3/2)B.$ Set $v:=(u-u_{2B})h.$ It is easily verified that $\widetilde g\in D(v),$ where
$$
\widetilde g:= {2\over r}|u-u_{2B}|\chi_{2B}+g\chi_{2B}.
$$
By $(q,q)$-Poincare, i.e. Proposition \ref{pr:2.4} combined with the H\"older-inequality, we have that $\| v\|^q_{L^q(\real^n)}\lesssim \int_{3B}g^q\, dx.$
 Now $v/u_{2B}$ is admissible for the definition of capacity, and we
may apply Theorem \ref{th:main} and the definition of capacity to
estimate
$$
cr^{n-q}\leq \capa_{H^{1,q}}(K;2B)
\lesssim (u_{2B})^{-q}\int_{\real^n}|\widetilde g|^q\, dx\lesssim (u_{2B})^{-q}\int_{3B}|g|^q\, dx. 
$$ 
This yields an upper bound for $u_{2B}$, which in combination with \refeq{eq:4.1011} 
finishes the proof.
\end{proof}

Let $p\in (0,1].$ We call  the complement of a domain  $\Omega\subset\real^n$ 
 uniformly $p$-fat if   there is a constant $c>0$
such that for each $x\in \Omega^c$ 
and $r>0$ one has 
\beqla{eq:4.1021} &&
\capa_{H^{1,p}}( \Omega^C\cap B(x,r);B(x,2r))\geq cr^{n-p}.
\eeq
Similarly, the complement $\Omega^c$ is thick in the sense of $s$- Hausdorff content (here $s>0$ is arbitrary) if there is $c>0$  such that for each $x\in \Omega^c$ 
and $r>0$ one has 
\beqla{eq:4.102} &&
{\mathcal H}_\infty^s(\Omega^c\cap B(x,r))\geq cr^{s}.
\eeq
Here for any subset $E\subset\real^n$ the $s$-Hausdorff content 
of $E$ is defined as ${\mathcal H}_\infty^s(E):=c(n,q)\inf\sum_{k=1}^\infty r_k^s$,
where the infimum is taken over all denumerable coverings of $E$
by balls with radii $r_k$, $k\in\nanu .$

\begin{theorem}\label{th:4.2} Let $p\in (n/(n+1),1]$. 
{\bf (i)}\quad 
Assume that $\Omega\subset \real^n$ is a bounded domain
with uniformly $q$-fat complement  for some $q<p$.
Then  the domain $\Omega $ carries the $p$-Hardy inequality {\rm \refeq{eq:4.40}}.

\smallskip

\noindent{\bf (ii)}\quad Assume that the complement of $\Omega$ is thick in the sense of $(n-q)$-Hausdorff content
with some $q<p$. 
Then  the domain $\Omega $ carries the $p$-Hardy inequality {\rm \refeq{eq:4.40}}.
\end{theorem}
\begin{proof}\quad (i)\quad
Let $u\in C_0^\infty (\Omega )$ and $g\in D(u).$
Fix  $x\in\Omega$. Denote $r=d(x,\partial\Omega )$ and
choose  $z\in\partial\Omega$ with $d(x,z)=r.$ We
may apply Lemma \ref{le:4.1} with the choice $B=B(z,r).$ The outcome is
$$
u(x)/d(x,\partial\Omega )\leq c(M(g^q)(x))^{1/q}.
$$
As $g^q\in L^{p/q},$ where $p/q>1$ we obtain the desired
conclusion by the Hardy-Littlewood maximal theorem.

\smallskip

(ii)\quad  It is clearly enough to prove the following: given
$n/(n+1)<q<p\leq 1$ and $C>0$ there is a constant $c=c(n,p,q)>0$
such that the condition $H^{n-q}_\infty (E)\geq Cr^{n-q}$
 implies that $\capa_{H^{1,p}}(E;2B)\geq cr^{n-p}.$
Here $B\subset\real^n$ is a ball and $E\subset \overline{B}$ is
an arbitrary compact subset. Actually, since both the capacity and the Hausdorff content scale in the appropriete manner,
it is enough to consider the case $r=1$. 

In order to prove the  claim  we thus assume that
$B$ is a ball
of radius one and $E\subset \overline{B}$ is
compact with $H^{n-q}_\infty (E)\geq c_0.$ Suppose that $u\in C_0^\infty (2B)$ satisfies $u\geq 1$ on $E$ and $g\in D(u).$ Fix an arbitrary
$x\in E$. Choose a ball $B_0\subset 3B$
with radius $1$ such that $B_0\cap 2B=\emptyset.$ Moreover,
for $k\geq 1$ denote $B_k=2^{1-k}B.$
By using again the Poincare inequality to compare mean values in the
chosen balls  we deduce that
\beqla{eq:4.120}
1=|u(x)-u_{B_0}|\leq c_1\left((\mint_{3B} g^q\,dx)^{1/q}+2^{-k}\sum_{k=1}^\infty(\mint_{B_k} g^q\,dx)^{1/q}\right).
\eeq
Since $\sum_{k=0}^\infty 2^{(-1+q/p)k}=c_2<\infty,$ we deduce that
necessarily either $(\smint_{3B} g^q\,dx)^{1/q}\geq 1/(c_1c_2)$
or there is an index $k\geq 1$ such that $(\smint_{B_k} g^q\,dx)^{1/q}\geq 2^{kq/p} /(c_1c_2).$ In both cases, by the
H\"older inequality, there is
a constant $c_3=c_3(p,q,n)$ such that  one of these balls,
call it $\widetilde B$, has  radius $\widetilde r\leq 3$, is centered at $x$  and satisfies
\beqla{eq:4.121}
\int_{\widetilde B}g^p\, dx\geq c_3{\widetilde r}^{n-q}.
\eeq

We select a ball $\widetilde B$ as above for each $x\in E.$ By a standard covering argument we may select a disjoint subcollection ${\mathcal B}$ such that the
corresponding 5-times inflated balls cover $E.$ By summing
the inequality \refeq{eq:4.121} over this 
collection we obtain
\beqla{eq:4.122}
\int_{\real^n}g^p\, dx\geq \sum_{\widetilde B\in{\mathcal B}}\int_{\widetilde B}g^p\,dx
\geq c_35^{q-n}H^{n-q}_\infty (E)\geq c_35^{q-n}c_0.
\eeq
This  concludes the proof of part (ii) of the Theorem.
 \end{proof}

\begin{remark}\label{re:4.10}
The proof of  part (ii)  is insipired by \cite{HajlaszKoskela2}. In connection with the above result we 
point out the interesting paper \cite{Orobitg} which connects
the $p$-Hardy Sobolev capacity directly to the $(n-p)$-Hausdorff content. It should be also noted that in the case $p>1$ one may
show the sufficiency of uniform $p$-fatness of $\Omega^c$ by employing the self-improving property
of uniform $p$-fatness, due to J. Lewis \cite{Lewis}, which approach appears not to work for $p\leq 1.$
We plan to return to Hardy's inequality in more general domains and consider at the same time generalizations of the results of the present paper to non-Euclidean situations.
\end{remark}


\begin{thebibliography}{kkk}

\bibitem{AucherEmmanuelTchamitchian}
P. Auscher, E. Russ, and P. Tchamitchian: {\it  
Hardy Sobolev spaces on strongly Lipschitz domains of $\mathbf R\sp n$,}
J. Funct. Anal. 218 
 (2005), no. 1, 54--109.

\bibitem{ChangKrantzStein} D.-C. Chang, S.G. Krantz, and E.M. Stein: {\it
 $H\sp p$ theory on a smooth domain in $R\sp N$ and elliptic boundary value problems,}
J. Funct. Anal. 114 (1993), 286--347.

\bibitem{ChangDafniStein} D.-C. Chang, G. Dafni, and E.M. Stein: {\it
Hardy spaces, BMO, and boundary value problems for the Laplacian on a smooth domain in $\bold R\sp n$,}
Trans. Amer. Math. Soc. 351 (1999), 1605--1661. 





\bibitem{Cho} Y.-K. Cho:  {\it Strichartz conjecture on Hardy-Sobolev spaces,} Colloq. Math. 103 (2005), 99-114.

\bibitem{Coifmann} R. R. Coifman: {\it A real variable characterization of $H^p$,} Studia Math. 51 (1974), 269--274. 

\bibitem{DeVoreSharpley}
R. DeVore and R. C. Sharpley:   {\it
Maximal functions measuring smoothness.} Mem. Amer. Math. Soc. 47 (1984), no. 293.


\bibitem{FeffermanStein} C. Fefferman and E. M. Stein: {\it $H^p$ spaces of several variables,} Acta Math. 129 (1972), 137--193. 

\bibitem{GaSeJi} A. Gatto, C. Segovia, and J. Jiménez: {\it On the solution of the equation $\Delta \sp{m}F=f$ for $f\in H\sp{p}$}. Conference on harmonic analysis in honor of Antoni Zygmund, Vol. I, II (Chicago, Ill., 1981),  394--415, 1983.

\bibitem{Goldberg} D. Goldberg: {\it A local version of Hardy spaces,} Duke Math. J. 46 (1979), 27--42.

\bibitem{Hajlasz1} P. Hajlasz: {\it  Sobolev spaces on an arbitrary metric space,}  Potential Anal.  5  (1996),  403--415

\bibitem{Hajlasz3} P. Hajlasz:  {\it Pointwise Hardy inequalities,} Proc. Amer. Math. Soc. 127 (1999), 417--423. 

\bibitem{Hajlasz2}  P. Hajlasz: {\it Sobolev spaces on metric-measure spaces.}  Heat kernels and analysis on manifolds, graphs, and metric spaces (Paris, 2002),  173--218, Contemp. Math., 338, Amer. Math. Soc., Providence, RI, 2003. 

\bibitem{HajlaszKoskela2}  P. Hajlasz and P.  Koskela: {\it Sobolev meets Poincar\' e,} C. R. Acad. Sci. Paris Sér. I Math. 320 (1995), no. 10, 1211--1215.

\bibitem{HajlaszKoskela}  P. Hajlasz and P.  Koskela: {\it Sobolev met Poincar\' e}.  Mem. Amer. Math. Soc.  145  (2000).
 
\bibitem{HaKoTu}  P. Hajlasz, P.  Koskela, and Heli Tuominen: {\it Measure density and extendability of Sobolev
functions}, manuscript. 

\bibitem{Janson} S. Janson: {\it On functions with derivatives in $H\sp 1$.} IN:  Harmonic analysis and partial differential equations (El Escorial, 1987),  193--201. Springer Lecture Notes in Math., 1384, 1989.

\bibitem{JoSjWa}
A. Jonsson, P. Sj\"ogren, and H. Wallin:
{\it Hardy and Lipschotz spaces on subdomains of $R\sp n$,}
Studia Math 80  (1984),  2, 141--166.


\bibitem{KiMa} J. Kinnunen and O. Martio: {\it Hardy's inequalities for Sobolev functions,} Math. Res. Lett. 4 (1997),  489--500.




\bibitem{Lewis} 
 J. Lewis: {\it Uniformly fat sets.}  Trans. Amer. Math. Soc.  308  (1988),  177--196.

\bibitem{Mi1}
A. Miyachi: {\it  $H\sp p$ spaces over open subsets of $R\sp n$,} Studia Math 95  (1990), 205--228.

\bibitem{Mi2}
A. Miyachi:  {\it Extension theorems for real variable Hardy and Hardy-Sobolev spaces.} Harmonic analysis (Sendai 1990),  170--182, ICM-90 Satell. Conf. Proc., Springer, Tokyo 1991. 

\bibitem{Mi3}
A. Miyachi:  {\it Hardy-Sobolev spaces and maximal functions,} J. Math. Soc. Japan 42 (1990), no. 1, 73--90.

\bibitem{Orobitg} J. Orobitg: {\it  Spectral synthesis in spaces of functions with derivatives in $H\sp 1$},  Harmonic analysis and partial differential equations (El Escorial, 1987),  202--206. Springer Lecture Notes in Math., 1384, 1989.

\bibitem{Peetre} J. Peetre: {\it New thoughts on Besov spaces.} Duke Univ. Math. Ser, Duke University 1976.

\bibitem{Stein} E. M. Stein: {\it Harmonic Analysis.} Princeton University Press 1993.

\bibitem{Strichartz} R. Strichartz:  {\it $H^p$ Sobolev spaces,} Colloq. Math. 60/61 (1990), 129--139.

\bibitem{Torchinsky} A. Torchinsky: {\it Real-variable methods in Harmonic analysis.} Academic Press 1986. 

\bibitem{Triebel} H. Triebel: {\it Theory of function spaces.} Monographs in Math. 78. Birkh\"auser 1983.






\end{thebibliography}
\end{document}